# THE SIMPLICITY OF COMPLEXITY, A STORY OF MATHEMATICAL OUTREACH

**Hugo Parlier (University of Fribourg) & Bruno Teheux (University of Luxembourg)**

Mathematics is notoriously hard to communicate. Communicating research level mathematics is near impossible, because the staggering pace of mathematical discovery leaves the general audience by and large in the dark. The notion of general audience doesn't really exist anyways, because scientific outreach needs to be tailored to specific audiences. In particular, cultural divides, including language and social barriers, make good practice more difficult to share. While all of these difficulties are often given as warnings for eager mathematicians wanting to share the beauty of their passion to a wider audience, despite being based in a mix of common sense and experience, maybe they're not as unsurmountable as they might seem.

This is an account of our experience, in part lead by chance and opportunity, in creating activities and bringing them to life in unexpected venues, including game festivals and two world expos.

**The beginnings**

We began working together in 2017, as colleagues in Luxembourg, via what one might call a "standard channel". We proposed a mathematical exhibit "Unpuzzling mathematics" for the Luxembourg Science Festival. Our proposal was initially refused because, from the refusal letter, "we had to make a choice between two mathematical workshops". We expressed our disappointment, and the organizers graciously offered us a free spot when there was a cancelation. The fact was that proposals from mathematics were rare, and it wasn't clear to the organizers how one could make mathematics attractive in a festival with cool lasers, robots and volcanoes to compete with.

That first festival was an eye-opener for us. Determined to make our exhibit a success, we aimed to create a space with a lounge type atmosphere, where visitors could enjoy mathematics without formalism, and so we presented a collection of games and activities, some that we had created and others that we bought, that had mathematics in the background. The games were an excuse to talk about optimization, graph theory and topology, but most of all to give visitors the experience of mathematical exploration. The exhibit was hugely successful, ensuring our spot among the science outreach community in Luxembourg who welcomed and encouraged us to do more.

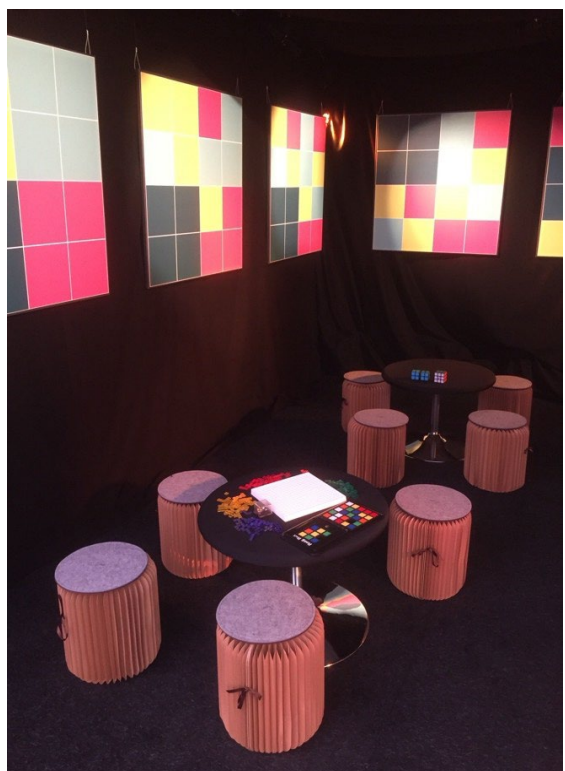

*Illustration 1. Our first science festival setup (2017)*

Emboldened by this first experience, we systematically created exhibits for science festivals. Without going into too much detail, we began to systemically create activities, always rooted in research and always tailor-made for the space we had imagined. While at first, we also used some of the wonderful games and activities that we found elsewhere, little by little all of our activities became homemade. For one thing, it is exhilarating to create new material but, more importantly, when you do, you can mould it to be able to best share whatever message or experience you are trying to get across. You can also design it to fit your atmosphere. Of course we are not interior designers, nor designers full stop, but making

Hugo Parlier (University of Fribourg) & Bruno Teheux (University of Luxembourg)

and creating spaces that we enjoy makes it easier to be enthusiastic when trying to spike other people's interest.

Retrospectively, "breaking through" to an audience that attends science festivals isn't particularly difficult, even if you are competing with volcanoes. The audience is naturally curious and eager to engage. Our lounge atmosphere, in the hustle and bustle of a festival full of visitors, was a welcome haven of tranquillity. Our activities were never labelled as mathematics, and by proposing games we appealed to visitors' playful nature. The trick was to make sure visitors were surprised by what they saw, that they left thinking differently about mathematics and, above all, with a smile on their face. Having large numbers of visitors is radical for testing activities. You see what works well, and, importantly, what doesn't and what requires more guidance. We always tried new activities alongside others that we knew worked. What quickly emerged from all of this, were guidelines that we had been implicitly following, but that moulded our way of creating and adapting new material.

**Principles**

Here are just a few of the principles that we follow (sometimes):

**Atmosphere**: As mentioned above, the first thing was that we always try and create an inviting atmosphere. As anything related to mathematics requires some amount of thinking, a relatively calm atmosphere works nicely. By creating more an artsy vibe, we also suggest something along the lines of "mathematics is an art" and, without getting into the debate of how true this is, certainly creativity is a vital ingredient in mathematical exploration. This also, conveniently, challenges what visitors might expect from a science that is sometimes perceived as dry, cold and unforgiving.

**Simplicity:** when faced with visitors who either want to go see the volcano or who have been afraid of mathematics since 6$^{th}$ grade – there are lots of these, even in science festivals –, you can either let them go and focus on those who are bolder, or try to hook them somehow. But you have very little time to do so. We came up with the 10 second rule: the rules of a game, or "how to" process of an activity, have to be deliverable in under 10 seconds. Mathematics is not a spectator sport: to enjoy them, you have to engage.

Simple and natural rules allow visitors to think "well, what can I lose by trying". By building from simple rules, the contrast makes the emergence of complexity even more stunning.

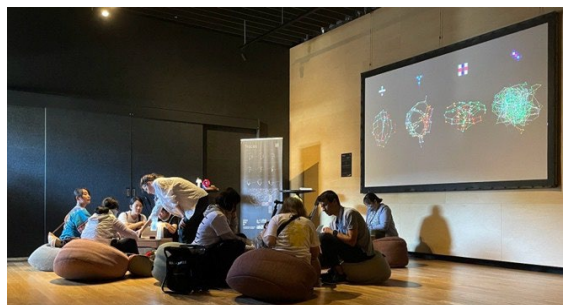

*Illustration 2. Our mathematical lounge in Osaka*

**Mediation**: The word "mediation" is less used in this context in English, in favour of "outreach", but it appeals to a vital part of the process. No matter how simple you make an activity, visitors regularly require encouragement. Sometimes it is a hint, sometimes a short conversation about "how is this math?", and often just a "you can do this". While a kid plays, ask a grandparent if they want to give it a try. In our experience, the kids sometimes have to drag the parents out to see the volcanoes.

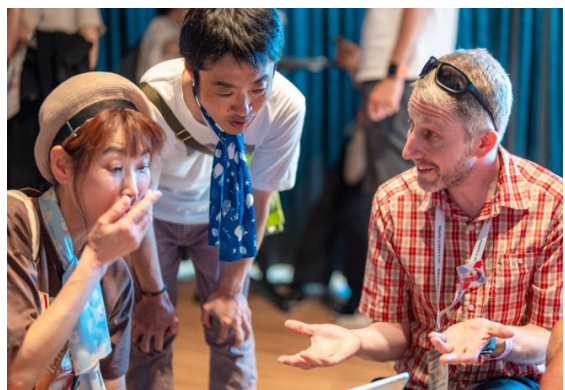

*Illustration 3. Mediation is key and leads to shared joy*

**Level up but go deep:** In particular, start easy. If you don't have the pleasure of success pretty quickly, this will only reinforce deeper fears about "math is not for me". It is remarkable how much time and effort visitors are willing to invest into more complicated ventures once they are emboldened by success and curiosity. How far or how deep visitors are willing to go varies of course. By creating activities that can be made

Hugo Parlier (University of Fribourg) & Bruno Teheux (University of Luxembourg)

progressively difficult, with natural "stopping points", visitors can go away with the satisfaction of having done something they had never done before, but also with the knowledge that there is more to explore.

Explicitly or implicitly hinting at research level math is a real bonus. If the activity provokes questions from visitors for which the only answer is "We wish we knew how to answer that: that's an open problem!", even better.

As with any set of principles, we are also happy to ignore them on occasion, especially for activities or visualizations that we reserve for the more tenacious visitors, those eager to know more about some of the deeper mathematics behind. The trick is not to over-explain. Not everyone who asks really wants to know all of the details. Spike curiosity and guide visitors willing or compelled to dig deeper, but don't force information upon them. As educators, it is sometimes hard not to share what we know. After all, in some sense, that's what we are trained to do. Empathy is more important than divulging knowledge in these contexts (and not only).

**New venues**

As our activities gained some traction, we started to get opportunities outside of science festivals and those aimed at schools. In particular, we were invited to participate in a game festival (the *Game On* festival in Luxembourg city). Again, a new challenge awaited. This time we weren't going to be engaging with a public who signed up for science. And having games in a science festival is original, but hardly in this context. Luckily, some of our guiding principles came to the rescue. There are many extraordinary games in the world, and the games presented at the festival were as one might expect, very cool. However, many good board games need some time to get into. Our 10 second rule works really well in this context: you can get visitors hooked very quickly. Plus, because we presented tailor-made games, most visitors have never seen anything like them.

And while, of course, there are science buffs who are also gaming enthusiasts, the audience as a whole was quite different in composition, with visitors and families in particular not necessarily having an academic background.

Right before our invitation to a game festival, through an extraordinary series of circumstances, we had been invited to exhibit some of our activities at the World Expo in Dubai in 2021/2022 on the Luxembourg Pavilion[1]. We had been eyeing the pavilion for a while, because we had heard that it had been designed in the shape of a Möbius band[2] which we thought was pretty cool (and it really was).

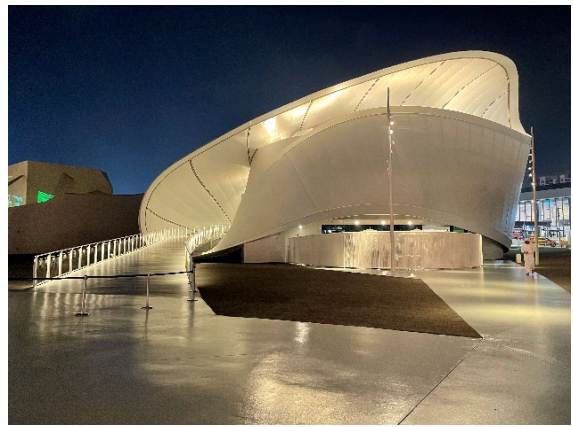

***Illustration 4.*** *Luxembourg Pavilion at Expo 2020 Dubai by Metaform Architects*

Via colleagues at Sorbonne Abu Dhabi, we met with the Luxembourg ambassador to the Emirates and associates, who were enthusiastic and put us in contact with the pavilion who had very nice space for temporary exhibits and invited us for a "learning themed" week. Once again, we were lucky in many regards. The Luxembourg Pavilion was very successful, in particular because they had a really fun slide, so there were lots of visitors and they could easily visit our exhibit along their route. Our motto became "Come for the slide, stay for the math", and, with something in the order of a thousand active daily visitors, the exhibit was a success. So much so that the pavilion asked if we would be willing to come back for a second residency, this time for a longer period. We did, and this time with a double purpose.

---

[1] https://gouvernement.lu/en/dossiers.gouv2024_mindigital+en+dossiers+2022+vr-dubai-2020.html.

[2] A Möbius band is obtained by pasting two opposite sides of a rectangle with a twist. It is an important example of a non-orientable surface and, in particular, it only has one side!

Hugo Parlier (University of Fribourg) & Bruno Teheux (University of Luxembourg)

In addition to games, we presented early models of the "Life Lines" experience that will be described below.

*Playing or creating. Competing or creating something collaborative. Solo or collective efforts.* Again, by means of contrast and unexpected results, this invited visitors to get a glimpse of the diversity that can lie within mathematics, or at least to view it through very different lenses. Again, the magic worked and with 20000 drawings collected, 50000 puzzles solved, the numbers were there to prove it.

Between the Expo, the game festival, and science festivals, we were convinced that we could reach an audience of all ages, and with a large diversity of backgrounds. Keeping things very simple, visual and self-explanatory makes for universal appeal. Just by the sheer number of visitors, we had learned by experience when to intervene, how much encouragement was necessary, and when to give the visitors space. We had also learned how to pick up on certain warning signs. Our main goal was not to channel people's math insecurities or to please parents eager to turn their kids into math wizzes.

When we were again contacted for the next universal expo, this time in Osaka, we were ready. In addition to a two-week residency on the Luxembourg Pavilion, we contacted the Swiss Pavilion and the Belgium Pavilion, both who were willing to host temporary events. Again, new challenges came up. Most participants were Japanese, and many spoke little or no English. While we didn't use a lot of words to explain our activities, we did need to say something. We were also warned that it was difficult to engage with participants, that the cultural barrier wasn't easy to breach, and that without Japanese translators we would get nowhere. We thus arranged for there to always be at least one person on site who spoke Japanese, but they weren't trained as mathematicians.

We were faced with the double difficulty of getting sophisticated mathematics, and without the possibility of detailed verbal explanations. Again, our principles came to the rescue. The rules of the games and activities are sufficiently intuitive that they can be demonstrated without words. And by sharing something new, surprising and which brings joy, allows to create an emotional bond that transcends cultural barriers. The universality of the pleasure of artistic creation or the satisfaction of solving a tricky puzzle is shared by visitor and mediator alike. We were also helped by the love of the Japanese for games and challenges. Not only would they engage, but every day we would have visitors, sometimes as young as 6 and sometimes well beyond retirement age, who would spend hours diligently trying to solve all the puzzles.

**So, what did visitors do?**

To give a flavour of the activities we presented, here are the activities we presented in Osaka.

Quadratis and Revolution

Quadratis[3] is a puzzle type game where you are presented with a pattern of square tiles, which are then shuffled, and your goal is to restore the original pattern by sliding them around. The shuffling and sliding moves are dictated by the way the squares are pasted together, so in fact you are playing on a surface with "topology", that is on a surface where sides are connected in possibly surprising ways. The puzzles range from downright easy to near impossible, with more or less everything in between. As an added advantage, visitors can download a version of the game for free.

Revolution[4] is also a reconfiguration puzzle, but this time played on a physical board. Here "moves" come from rotating a wheel to move rings around, and your goal is to find a prescribed sequence of ring patterns. As the game progresses, puzzles become increasingly challenging.

---

[3] Quadratis was co-created by H. Parlier and P. Turner with development by M. Gutierrez and R. Juarez and is available on http://quadratis.app.

[4] Revolution was created and designed by H. Parlier and B. Teheux.

Hugo Parlier (University of Fribourg) & Bruno Teheux (University of Luxembourg)

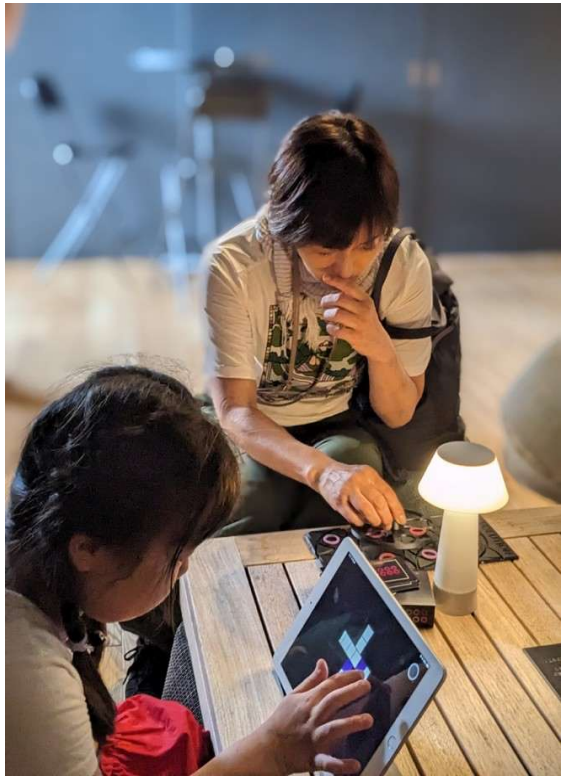

*Illustration 5. Playing Revolution and Quadratis in Osaka*

Behind both games are notions from topology, geometry and graph theory. Most questions one could ask about these games are largely unexplored: how many moves do I need to go from here to here? Is it always possible to relate two configurations? If so, how fast would an optimal algorithm do it? Can you write down an explicit algorithm to solve the puzzle?

These are all research level questions. In particular, some of these boil down to deep questions about moduli spaces[5], vast spaces which encode how shapes relate to one another, and hot research topics in geometry. Other questions are directly related to graph theory and computational geometry. Graphs – or networks of possible patterns – can be visualized tracking moves made by players.

### Life Lines

In Life Lines[6], visitors are invited to use their finger to draw a curve of any shape they please with the only constraint being that it goes from left to right. The resulting drawing is then mysteriously coloured. You can then add more colours or rearrange them, but you quickly find out that your control is limited.

There seem to be laws dictating how zones are coloured. The drawings are collected and then shared with other visitors through interactive art pieces. Visitors can also explore the space of all drawings made by other visitors, and create morphings to see how two drawings relate.

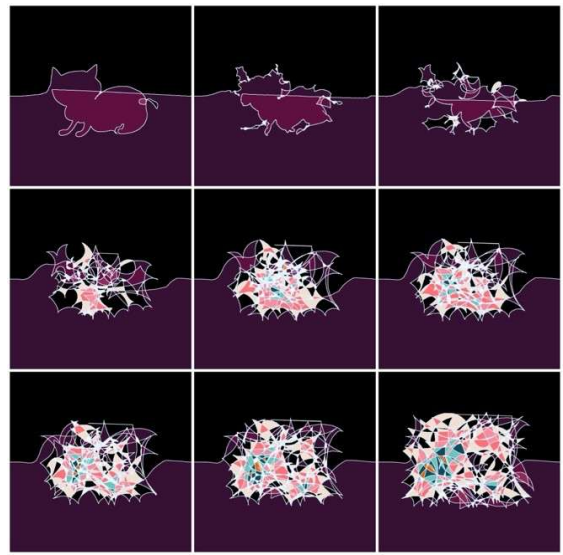

*Illustration 6. A morphing between figurative and abstract creations*

How *are* the picture coloured? We use another notion from topology, the so-called winding number, to track how many times curves "wrap" around the zones delimited by the drawing, taking into account the orientation of the drawing. This leads to a type of natural harmony among colours.

---

[5] Moduli spaces are vast structured landscapes of all possible geometric forms. One of the most important experts in moduli space theory, and its relationship to curves on surfaces and hyperbolic geometry, was Maryam Mirzakhani (1977–2017), the first female Fields medalist. For her, mathematics were inherently an artistic pursuit: "The beauty of mathematics only shows itself to more patient followers."

[6] Life Lines, created by H. Parlier and B. Teheux, was first showcased in its current form at the *Shapes* exhibit at EPFL Pavilions in January 2025. See https://epfl-pavilions.ch/en/exhibitions/shapes.

And again, this activity is directly related to hot research topics. For one thing, from a philosophical point of view at least, by creating a drawing, you are sampling the space of all possible drawings, which again is a type of moduli space. But how do you measure proximity between drawings? We use what is commonly called the "dog-leash distance" (and formally called Fréchet distance) to measure how close two drawings are. We use recent efficient algorithms to do this, but there would be many other ways to measure these distances. With tens of thousands of drawings, we now have an extraordinary data base of collective creativity to explore using modern geometry.

**What's next? Why should we be doing this in the first place?**

There is a very different point of view when doing mathematical outreach which is to insist on how important mathematics is, as a language for science and a tool for understanding the world.

Our point of view is hardly opposite to that, but it is definitely orthogonal. The activities we present are, hopefully, intrinsically interesting. It might be effective to play basketball to keep fit, but you can also play it just to have fun. Having fun doesn't prevent you from improving your cardio. Reaching out is also about reaching a new public, people who aren't aware of how interesting mathematics can be. Not every visitor needs to be turned into a mathematician, but one thing that would be great, is that every visitor leaves with the idea that mathematics is a thriving science with lots to be explored.

Department of Mathematics, University of Fribourg, Chemin du Musée 23, 1700 Fribourg, Switzerland.

E-mail: hugo.parlier@unifr.ch

Website: https://homeweb.unifr.ch/parlierh/pub/

ORCID: 0000-0001-5618-509X

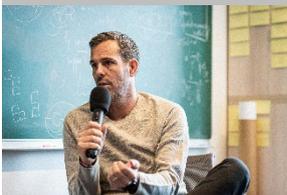

**Hugo Parlier**, Full Professor at the University of Fribourg since 2025, is a mathematician interested in the study of shapes and their deformations. Many of his results are related to the study of curves on surfaces and often have a combinatorial and visual flavour.

In addition, he is passionate about sharing the ever-evolving nature of mathematical research with as wide an audience possible.

Recently, his efforts and fruitful collaborations, often in collaboration with Bruno Teheux, led to activities involving puzzles and collaborative art activities being showcased in science museums, at the EPFL Pavilions, and at the World Expos in 2022 and 2025.

*Photo credentials: Mike Zenari*

Department of Mathematics, University of Luxembourg, Maison du Nombre, Avenue de la Fonte 6, 4364 ESCH, Luxembourg.

E-mail: Bruno.teheux@uni.lu

Website: https://math.uni.lu/teheux/

ORCID: 0000-0002-3007-3089

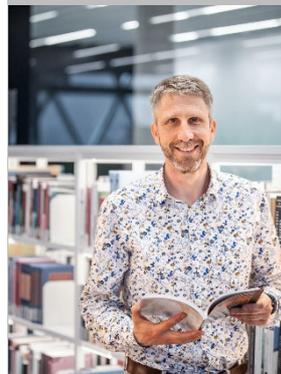

**Bruno Teheux**, Assistant professor at the University of Luxembourg, is a mathematician active in mathematical logic. He investigates the semantics of non-classical logics through the lens of algebra and topology.

With his outreach activities, he contributes to creating and mediating at science festivals, games events, cultural festivals, and world expositions, where he loves updating audiences' perspectives on mathematics. His goal is to combine elements such as games or art with mathematics to create a lasting emotion.

*Photo credentials: Remi Grizard*